\documentclass[a4,12pt]{article}
\usepackage{amsfonts}
\usepackage{bbding}
\usepackage{geometry}
\usepackage{multicol}
\usepackage{graphicx}
\usepackage{subfigure}
\usepackage{amsmath}
\usepackage{amsmath,bm}
\usepackage{mathrsfs}
\usepackage{amssymb}

\usepackage{amsfonts}
\usepackage{bbding}
\newtheorem{thm}{Theorem}[section]
\newtheorem{lem}[thm]{Lemma}
\newtheorem{rmk}[thm]{Remark}
\newtheorem{defi}[thm]{Definition}
\newtheorem{ques}[thm]{Question}
\newtheorem{eg}[thm]{Example}
\def\deg{\hbox{\rm{deg\,}}}

\def\supp{\hbox{\rm{supp}}}

\renewcommand{\vec}[1]{\bm{#1}}
\def\tr{{\rm{tr}}}
\newcommand{\wt}{\hbox{\rm{wt}}}

\def\proof{\noindent{\textbf{Proof}}\quad}
\newcommand{\qedd}{\hspace*{\fill}$\Box$\medskip}   

\begin{document}

\title{\LARGE{\textbf{$\mathcal{P}\mathcal{S}$ bent functions constructed from finite pre-quasifield
spreads}}}
\author{\large{Baofeng Wu}}
\date{\normalsize State Key Laboratory of Information Security, Institute of Information Engineering,\\  Chinese Academy of Sciences,
Beijing 100093, China} \maketitle

\noindent\textbf{Abstract}\quad Bent functions are of great
importance in both mathematics and information science. The
$\mathcal{P}\mathcal{S}$ class of bent functions was introduced by
Dillon in 1974, but functions belonging to this class that can be
explicitly represented are only the
$\mathcal{P}\mathcal{S}_{\text{ap}}$ functions, which were also
constructed by Dillon after his introduction of the
$\mathcal{P}\mathcal{S}$ class. In this paper, a technique of using
finite pre-quasifield spread from finite geometry to construct
$\mathcal{P}\mathcal{S}$ bent functions is proposed. The constructed
functions are in similar styles with the
$\mathcal{P}\mathcal{S}_{\text{ap}}$ functions. To explicitly
represent them in bivariate forms, the main task is to compute
compositional inverses of certain parametric permutation polynomials
over finite fields of characteristic 2. Concentrated on the
Dempwolff-M\"uller pre-quasifield, the Knuth pre-semifield and the
Kantor pre-semifield, three new subclasses of the
$\mathcal{P}\mathcal{S}$ class are obtained. They are the only
sub-classes that can be explicitly constructed more than 30 years
after the $\mathcal{P}\mathcal{S}_{\text{ap}}$ subclass was
introduced.\medskip

\noindent\textbf{Keywords}\quad Bent function;   partial spread;
finite prequasifield;   permutation polynomial;   compositional
inverse

\section{Introduction}\label{secintro}

Bent functions were introduced by Rothaus in 1976 {\cite{rothaus}}.
Roughly speaking, a bent function is a Boolean function whose Walsh
transform has a constant magnitude. Since they were introduced, bent
functions have attracted a lot of attention both in theory and
application as they are of interest in various aspects of
mathematics and information science. In combinatorics, they are
related to difference sets in elementary Abelian 2-groups and
Hadamard matrices. In cryptography, they are functions attaining
maximal nonlinearity, which have best resistance to the powerful
linear attack when used as primitives in designing such key stream
generators as filter generator and combiner generator in stream
cipher. In the coding context, they are equivalent to deep holes of
a well-known class of code called Reed-Muller code. Therefore,
constructions of bent functions, especially large classes of bent
functions, are very meaningful and have been extensively studied.
However, it is far from simple to construct bent functions and there
are only a few constructions are known. We refer to
\cite{carlet2010} for a survey of this topic.

In his PhD thesis in 1974, Dillon proposed a construction of a very
large class of bent functions known as the partial spread class, or
$\mathcal{P}\mathcal{S}$ class for short {\cite{dillion1974}}. He
also studied a subclass of it named the
$\mathcal{P}\mathcal{S}_{\text{ap}}$ class, functions belonging to
which could be explicitly represented by bivariate polynomials over
finite fields. In fact, Dillon's construction of
$\mathcal{P}\mathcal{S}$ bent functions is more like a principle of
construction than an explicit construction. In other words, Dillon
just showed how to construct bent function along his rule, but did
not expressed what these functions are like, except the functions in
the $\mathcal{P}\mathcal{S}_{\text{ap}}$ class. Thus the following
question is natural.

\begin{ques}
How to explicitly represent other subclasses of the
$\mathcal{P}\mathcal{S}$ class except the
$\mathcal{P}\mathcal{S}_{\text{ap}}$ subclass by (univariate or
multivariate) polynomials over finite fields?
\end{ques}

In fact, it is not a simple matter to answer this question. To the
best of the author's knowledge, the
$\mathcal{P}\mathcal{S}_{\text{ap}}$ class is the only subclass of
the $\mathcal{P}\mathcal{S}$ class that can be explicitly
represented up to present. That means no progress in explicitly
constructing $\mathcal{P}\mathcal{S} $ bent functions has been made
in more than 30 years since the $\mathcal{P}\mathcal{S}$ and
$\mathcal{P}\mathcal{S}_{\text{ap}}$ classes were introduced.

In this paper, we devote to giving an answer to the above question.
By carefully examining the structure of partial spreads, we find
that a subclass of the $\mathcal{P}\mathcal{S}$ class can be derived
from an algebraic structure called pre-quasifield, which is an
important tool in studying finite geometry (geometry over finite
fields). To explicitly represent functions of a subclass, the main
difficulty is to compute compositional inverse of a parametric
permutation polynomial determined by the multiplication of the
pre-quasifield deriving this subclass. Generally speaking, to
compute compositional inverses of permutation polynomials over
finite fields is not easy, and there are rare classes of permutation
polynomials whose compositional inverses can be explicitly
determined. But for some known pre-quasifields, namely, the
Dempwolff-M\"uller pre-quasifield, the Knuth pre-semifield and the
Kantor pre-semifield, we can explicitly derive the compositional
inverses of the parametric permutation polynomials determined by
their multiplications. Hence we obtain three new subclasses of the
$\mathcal{P}\mathcal{S}$ class that can be explicitly represented by
bivariate polynomials over finite fields.

The rest of the paper is organized as follows. In Section
\ref{secpreli}, we recall some necessary backgrounds. In Section
\ref{secconst}, we propose the general technique to construct
$\mathcal{P}\mathcal{S}$ bent functions from pre-quasifield spreads.
As special examples, three new subclasses of the
$\mathcal{P}\mathcal{S}$ class are explicitly obtained  in Section
\ref{sec3class}. Concluding remarks are given in Section
\ref{secconclu}.

\section{Notations and preliminaries}\label{secpreli}

Let $n$ be a positive integer and $\mathbb{F}_{2^n}$ be  the finite
field with $2^n$ elements. For any $x\in\mathbb{F}_{2^n}$, we
sometimes distinguish $\frac{1}{x}$ with $x^{2^n-2}$. Denote by
``$\tr^n_1$" the trace function from $\mathbb{F}_{2^n}$ to
$\mathbb{F}_{2}$, i.e. $\tr_1^n(x)=\sum_{i=0}^{n-1}x^{2^i}$ for any
$x\in\mathbb{F}_{2^n}$.

\subsection{Boolean functions and their representations}

An $n$-variable Boolean function is a map from $\mathbb{F}_{2^n}$ to
$\mathbb{F}_{2}$. We denote by $\mathbb{B}_{n}$ the set of all
$n$-variable Boolean functions. The support of a Boolean function
$f\in\mathbb{B}_n$ is defined by
\[\supp(f)=\{x\in\mathbb{F}_{2^n}\mid f(x)=1\},\]
and the cardinality of $\supp(f)$ is often known as the Hamming
weight of $f$, denoted by $\wt(f)$. $f$ is called balanced if
$\wt(f)=2^{n-1}$.

By Lagrange interpolation, we can represent $f$ by a univariate
polynomial over $\mathbb{F}_{2^n}$ of the form
\[f(x)=\sum_{i=0}^{2^n-1}a_ix^i.\]
$\deg f=\max_{0 \leq i\leq 2^n-1}\{\wt(i)\mid a_i\neq0\}$ is called
the algebraic degree of $f$, where $\wt(i)$ represents the Hamming
weight, i.e. the number of 1's in the binary expansion, of the
integer $i$. In fact, representations of Boolean functions are more
flexible than what we can fully describe. For example, when $n=2m$
is an even integer, we can represent $f$ by  a bivariate polynomial
over $\mathbb{F}_{2^m}$ of the form
\[f(x,y)=\sum_{i=0}^{2^m-1}\sum_{j=0}^{2^m-1}a_{ij}x^iy^j\]
according to the isomorphism
$\mathbb{F}_{2^n}\cong\mathbb{F}_{2^m}\times\mathbb{F}_{2}$ between
$\mathbb{F}_2$-vector spaces. Under this bivariate representation,
the algebraic degree of $f$ is $\deg f=\max_{0 \leq i,j\leq
2^m-1}\{\wt(i)+\wt(j)\mid a_{ij}\neq0\}$ {\cite{carlet2010}}.

\subsection{Bent functions and Dillon's $\mathcal{P}\mathcal{S}$ class}

Let $f\in\mathbb{B}_n$. The Walsh transform of $f$ at
$w\in\mathbb{F}_{2^n}$ is defined by
\[\hat{f}(w)=\sum_{x\in\mathbb{F}_{2^n}}(-1)^{f(x)+\tr_1^n(wx)}.\]
When $n=2m$ and we use the bivariate representation of $f$, its
Walsh transform at $(u,v)\in\mathbb{F}_{2^m}\times\mathbb{F}_{2^m}$
can be expressed as
\[\hat{f}(u,v)=\sum_{(x,y)\in\mathbb{F}_{2^m}\times\mathbb{F}_{2^m}}(-1)^{f(x,y)+\tr_1^m(ux+vy)}.\]

\begin{defi}{{\cite{rothaus}}}\label{defbent}
A Boolean function $f\in\mathbb{B}_n$ is called a bent function if
$\hat{f}(w)=2^{n/2}$ for any $w\in\mathbb{F}_{2^n}$.
\end{defi}

It is implied by Definition \ref{defbent} that the number of
variables of  a bent function must be even.

Now we assume $n=2m$. In \cite{dillion1974} Dillon proposed a
principle of constructing bent functions using partial spreads from
finite geometry. Firstly we give the definition of spreads over an
$\mathbb{F}_2$-vector space.

\begin{defi}
Let $V$ be an $n$-dimensional vector space over $\mathbb{F}_2$.
Assuming $V_1,~\ldots,~V_{2^m+1}$ are all $m$-dimensional subspaces
of $V$ satisfying that
\[V=\bigcup_{i=1}^{2^m+1}V_i\]
and for any $1\leq i\neq j\leq 2^m+1$,
\[V_i\cap V_j=\{\vec0\}.\]Then $\mathfrak{S}=\{V_i\}_{i=1}^{2^m+1}$
is called a spread of dimension $m$ over $V$ and each element of
$\mathfrak{S}$ is called a component of it.
\end{defi}

\begin{eg}\label{spreadeg1}
For any $a\in\mathbb{F}_{2^m}$, let
\[V_a=\{(x,ax)\mid x\in\mathbb{F}_{2^m}\}.\]Besides, let
\[V_\infty=\{(0,x)\mid x\in\mathbb{F}_{2^m}\}.\]
It is direct to prove that $\mathfrak{S}=\{V_a\mid
a\in\mathbb{F}_{2^m}\cup\{\infty\}\}$ becomes an $m$-dimensional
spread of $\mathbb{F}_{2^m}\times\mathbb{F}_{2^m}$, which is an
$n$-dimensional vector space over $\mathbb{F}_2$.
\end{eg}

For a given $m$-dimensional spread
$\mathfrak{S}=\{V_i\}_{i=1}^{2^m+1}$ of $\mathbb{F}_{2^n}$ (viewed
as an $\mathbb{F}_2$-vector space), Dillon proposed the following
construction of $n$-variable Boolean functions: choose any $2^{m-1}$
components of $\mathfrak{S}$, say, $V_1,\ldots,V_{2^{m-1}}$ without
loss of generality, and set
\[f(x)=\sum_{i=1}^{2^{m-1}}\bm{1}_{V_i}\mod2,\]
where $\bm{1}_S(x)$ represents the character function on a set $S$,
i.e.
\[\bm{1}_S(x)=\left\{\begin{array}{ll}
1&\text{if}~x\in S,\\
0&\text{if}~x\not\in S.
\end{array}\right.\]
The function $f$ is named a $\mathcal {P}\mathcal {S}^-$ function by
Dillon. He similarly  defined a $\mathcal {P}\mathcal {S}^+$
function by choosing $2^{m-1}+1$ instead of $2^{m-1}$ components of
$\mathfrak{S}$ in the  construction  and called a function that was
$\mathcal {P}\mathcal {S}^-$ or $\mathcal {P}\mathcal {S}^+$ a
$\mathcal {P}\mathcal {S}$ function. It can be proved that all
$\mathcal {P}\mathcal {S}$ functions are bent. Obviously $f+1$ is a
$\mathcal {P}\mathcal {S}^+$ function for a $\mathcal {P}\mathcal
{S}^-$ function $f$, so we just need to talk about $\mathcal
{P}\mathcal {S}^-$ functions in the following.

It is clear that Dillon's construction of $\mathcal {P}\mathcal
{S}^-$ bent function is just a principle of construction since there
are various spreads over $\mathbb{F}_{2^n}$. If the spread is fixed
to be the $\mathfrak{S}$ in Example \ref{spreadeg1}, the
corresponding $\mathcal {P}\mathcal {S}^-$ functions can be
explicitly represented by  bivariate polynomials over
$\mathbb{F}_{2^m}$ of the form
\[f(x,y)=g\left(\frac{y}{x}\right),\]
where $g$  is an $m$-variable balanced Boolean function (represented
in univariate form) satisfying $g(0)=0$. All such functions form a
subclass of the $\mathcal {P}\mathcal {S}^-$ class named $\mathcal
{P}\mathcal {S}_{\text{ap}}$ class by Dillon. It is the only
subclass of the $\mathcal {P}\mathcal {S}^-$ class that can be
explicitly represented up to now.

\subsection{Finite pre-quasifield  and its spread}

\begin{defi}{{\cite{johnson2007}}}\label{defpqf}
A binary system $(Q,+,\diamond)$ is called a finite left
pre-quasifield, if $|Q|<\infty$ and the follow axioms holds:

\noindent(1) $(Q,+)$ is an Abelian group;

\noindent(2) $(Q^*,\diamond)$ is a quasigroup where
$Q^*=Q\backslash\{0\}$ (0 is the identity of $(Q,+)$), i.e. for any
$a\in Q^*$, the left multiplication operator
\begin{eqnarray*}
  L_a:~Q^* &\longrightarrow& Q^* \\
   x&\longmapsto&a\diamond x
\end{eqnarray*} and the right
multiplication operator \begin{eqnarray*}
  R_a:~Q^* &\longrightarrow& Q^* \\
   x&\longmapsto&x\diamond a
\end{eqnarray*}
are both bijective;

\noindent(3) For any $x,~y,~z\in Q$,
\[x\diamond(y+z)=x\diamond y+x\diamond z;\]

\noindent(4) For any $x\in Q$, $0\diamond x=0$.
\end{defi}

Similarly we can define finite right pre-quasifield by replacing
axiom (3) in Definition \ref{defpqf} by ``for any $x,~y,~z\in Q$,
$(y+z)\diamond x=y\diamond x+z\diamond x$" and axiom (4) by ``for
any $x\in Q$, $x\diamond 0=0$". Furthermore, $(Q,+,\diamond)$ is
called a pre-semifield if it is simultaneously  a left and right
pre-quasifield. In fact, for a right pre-quasifield $(Q,+,\star)$,
we can obtain a left pre-quasifield $(Q,+,\diamond)$ by defining
\[x\diamond y=y\star x,~\forall x,~y\in Q,\]
thus we need only to talk about left pre-quasifield in the following
and omit ``left" for simplicity. In addition, it can be proved that
the additive group of a pre-quasifield must be elementary Abelian,
so we can define a pre-quasifield by defining a multiplication in a
vector space over a finite field, whose characteristic is also
defined to be the characteristic of the corresponding
pre-quasifield.

Now we assume $(\mathbb{F}_{2^m},+,\diamond)$ is a pre-quasifield.
Set \[E_\infty=\{(0,x)\mid x\in\mathbb{F}_{2^m}\}\] and for any
$a\in\mathbb{F}_{2^m}$, set
\[E_a=\{(x,a\diamond x)\mid x\in\mathbb{F}_{2^m}\}.\]
According to the axioms in Definition \ref{defpqf}, it can be proved
that $\mathfrak{E}=\{E_a\mid a\in\mathbb{F}_{2^m}\cup\{\infty\}\}$
is a $m$-dimensional spread over
$\mathbb{F}_{2^m}\times\mathbb{F}_{2^m}$ {\cite{johnson2007}}. This
spread is often known as a pre-quasifield spread. Further results
shows that for any spread $\mathfrak{S}$ over
$\mathbb{F}_{2^m}\times\mathbb{F}_{2^m}$, there always exists a
pre-quasifield such that $\mathfrak{S}$ can be induced by it in the
above manner {\cite{johnson2007}}.

\subsection{Permutation polynomials over finite fields and their compositional inverses}

Let $\mathbb{F}_{q}$ be the finite field with $q$ elements where $q$
is a prime or a prime power. For any polynomial $F(x)$ over
$\mathbb{F}_{q}$, it can induce a map from $\mathbb{F}_{q}$ to
itself. $F(x)$ is called a permutation polynomial if the map induced
by it is bijective {\cite{lidl}}. It is clear that under the
operation of composition of polynomials and subsequent reduction
modulo $(x^q-x)$, the set of all permutation polynomials over
$\mathbb{F}_{q}$ forms a group which is isomorphic to $\mathcal
{S}_q$, the symmetric group on $q$ letters. Hence for any
permutation polynomial $F(x)$ over $\mathbb{F}_{q}$, there always
exists a unique polynomial $F^{-1}(x)$ such that
$$F(F^{-1}(x))\equiv F^{-1}(F(x))\equiv x \mod (x^q-x).$$ $F^{-1}$ is
called the compositional inverse of $F$ (or vice versa).

Permutation polynomials have many important applications  in
cryptography, coding theory and combinatorics. However, it is far
from easy to construct them. There are only a few classes of
permutation polynomials known. See \cite{lidl2, lidl3} for a survey
of this topic and \cite{akbary2, charpin, xdhou}, for example, for
some recent progresses. On the other hand, for a given class of
permutation polynomials, it is simple to verify whether the other
class of permutation polynomials represent their compositional
inverses, but it seems  very difficult  to explicitly derive their
compositional inverses. In fact, there are rare classes of
permutation polynomial whose compositional inverses can be
explicitly determined. We refer to \cite{coulter,wu1,wu2} for some
results on this topic. We involve the following results which will
be used in Section \ref{sec3class}.

\begin{thm}{{\cite{lidl4}}}\label{invdickson}
Let\[D_k(x)=\sum_{i=0}^{\lfloor\frac{k}{2}\rfloor}\frac{k}{k-i}{{k-i}\choose{i}}x^{k-2i}.\]
Then $D_k(x)$ is a permutation polynomial over $\mathbb{F}_q$ if and
only if $(k,q^2-1)=1$. Furthermore, if $k'$ satisfies
$kk'\equiv1\mod q^2-1$, then $D_k^{-1}(x)=D_{k'}(x)$.
\end{thm}

\begin{rmk}
The polynomial $D_k(x)$ in Theorem \ref{invdickson} is often called
a Dickson polynomial, which belongs to an important class of special
polynomials over finite fields.
\end{rmk}

\begin{thm}{\cite{wu3}}\label{invLP}
Let $m$ be odd and $a\in\mathbb{F}_{2^m}$ with $\tr(1/a)=1$. Then
\[L_a(x)=ax+a^2x^2+\tr(x)\]
is a permutation polynomial over $\mathbb{F}_{2^m}$. Furthermore,
the compositional inverse of $L_a(x)$ is
\[L_a^{-1}(x)=\frac{1}{a}C_a(x)+\frac{1}{a}\tr\left(\frac{x}{a}\right),\]
where $C_a(x)=\sum_{i=0}^{m-1}c_ix^{2^i}$ with coefficients given by
\begin{eqnarray*}
  c_0 &=& \frac{1}{a}\langle1,3,5,\ldots,m-3\rangle, \\
  c_i &=& \left\{\begin{array}{ll}
1+\frac{1}{a}\langle1,3,5,\ldots,i-2,\;i+1,i+3,\ldots,m-1\rangle&\text{if}~
i~ \text{is~odd},\\[.2cm]
\frac{1}{a}\langle0,2,4\ldots,i-2,\;i+1,i+3,\ldots,m-2\rangle&
\text{if}~ i~ \text{is~even},
  \end{array}
  \right.
\end{eqnarray*}
\noindent $1\leq i\leq m-1$.
\end{thm}

\begin{rmk}
The meaning of the symbol $\langle\cdot\rangle$ in Theorem
\ref{invLP} is as follows: for any $c\in\mathbb{F}_{2^m}$ and
$i_1,\ldots,i_s\in\mathbb{Z}/m\mathbb{Z}=\{0,1,\ldots,m-1\}$,
\[c\langle i_1,\ldots,i_s\rangle:=\sum_{j=1}^sc^{2^{i_j}}.\]
\end{rmk}

\section{Constructing $\mathcal
{P}\mathcal {S}$ bent functions form pre-quasifield
spreads}\label{secconst}

Let $(\mathbb{F}_{2^m},+,\diamond)$ be a pre-quasifield and
$\mathfrak{E}$ be its spread. For any
$a\in\mathbb{F}_{2^m}\cup\{\infty\}$, it is called the slope of the
component $E_a$. It is clear that each component of $\mathfrak{E}$
is fully determined by its slope.

Now we introduce a new operation on $(\mathbb{F}_{2^m},+,\diamond)$
named (left) division and denote it by
``$\diamond\frac{\cdot}{\cdot}$". For any
$a,~x\in\mathbb{F}_{2^m}^*$, $y\in\mathbb{F}_{2^m}$, we define
$\diamond\frac{y}{x}=a$ if $y=a\diamond x$. Besides, we fix
$\diamond\frac{y}{0}=0$. Then it is clear that $(x,y)\in E_a$ if and
only if $\diamond\frac{y}{x}=a$ for any $a\in\mathbb{F}_{2^m}$.

Recall the construction of $\mathcal {P}\mathcal {S}^-$ bent
functions. It is not difficult to see that the feature of a
$\mathcal {P}\mathcal {S}^-$ function  is that its limitation to all
nonzero elements of each component of the spread is a constant (0 or
1). Hence we obtain the following representation of $\mathcal
{P}\mathcal {S}^-$ bent functions.

\begin{thm}\label{mainthm}
Let $n=2m$ and $(\mathbb{F}_{2^m},+,\diamond)$ be a pre-quasifield.
Assum $g\in\mathbb{B}_{m}$ is a balanced  Boolean function with
$g(0)=0$. Then the function
\[f(x,y)=g\left(\diamond\frac{y}{x}\right)\]
is a $\mathcal {P}\mathcal {S}^-$ bent function.
\end{thm}

The proof of Theorem \ref{mainthm} is obvious and will be omitted.
Since every spread over $\mathbb{F}_{2^m}\times\mathbb{F}_{2^m}$ can
be induced by a pre-quasifield as mentioned before, we can represent
every $\mathcal {P}\mathcal {S}^-$ bent function in the form
indicated by Theorem \ref{mainthm} with the aid of the division
operation of the corresponding pre-quasifield. To explicitly
represent them, the rest problem is just to explicitly represent the
division operation by a polynomial over finite fields. For this
purpose, we need to examine the characterization of the
multiplication operation of a pre-quasifield.

Obviously for the pre-quasifield $(\mathbb{F}_{2^m},+,\diamond)$,
there exists a bivariate polynomial $F(x,y)$ over $\mathbb{F}_{2^m}$
such that $$x\diamond y=F(x,y).$$ We view $F(x,y)$ as a univariate
polynomial $F_y(x)$ by considering the variable $y$ to be a
parameter. Then $F_y(x)$ has the following property.

\begin{lem}\label{pqfmulti}
$F_0(x)=0$ and for any $y\in\mathbb{F}_{2^m}^*$, $F_y(x)$ is a
permutation polynomial over $\mathbb{F}_{2^m}$.
\end{lem}
\proof According to axiom (3) in Definition \ref{defpqf}, we have
\[F_0(x)=x\diamond 0=x\diamond(0+0)=x\diamond 0+x\diamond 0,\]
thus $F_0(x)=0$. Besides, for any $y\in\mathbb{F}_{2^m}^*$, if we
assume $F_y(x_1)=F_y(x_2)$, i.e.
\[x_1\diamond
y=x_2\diamond y,\] then  according to axiom (2) we get that
$x_1=x_2$. This is equivalent to say $f$ is a permutation
polynomial. \qedd

Now for any $a,~x\in\mathbb{F}_{2^m}^*$, $y\in\mathbb{F}_{2^m}$, if
$y=a\diamond x=F(a,x)=F_x(a)$, we can clearly get
\[\diamond\frac{y}{x}=a=F_x^{-1}(y).\]
That is to say if we want to explicitly represent the division
operation $\diamond\frac{y}{x}$ by a bivariate polynomial over
$\mathbb{F}_{2^m}$, we need only to compute the compositional
inverse of the parametric permutation polynomial $F_y(x)$ derived
from the multiplication operation of the pre-quasifield and reverse
the role of $x$ and $y$ in its formula afterwards.

\begin{eg}
Assume the multiplication of the pre-quasifield is just defined by
\[x\diamond y=F(x,y)=xy.\]
Then it is obvious that
\[\diamond\frac{y}{x}=F_x^{-1}(y)=\frac{y}{x}.\]
Hence the functions derived from this multiplication according to
Theorem \ref{mainthm} are just the
$\mathcal{P}\mathcal{S}_{\text{ap}}$ functions.
\end{eg}

To summarize, we have related the problem of constructing $\mathcal
{P}\mathcal {S}$ bent functions with the problem of computing
compositional inverses of special parametric permutation polynomials
over finite fields. Though it is challenging in general to determine
compositional inverse of a permutation polynomial, we can overcome
the difficulty if it is derived from the multiplication of certain
special pre-quasifield. In the next section we will focus on three
known pre-quasifields and explicitly represent the bent functions
derived from them.

\section{Three classes of $\mathcal
{P}\mathcal {S}$ bent functions constructed form pre-quasifield
spreads}\label{sec3class}

In this section, we talk about three well-known classes of
pre-quasifields of characteristic 2, namely, the Dempwolff-M\"uller
pre-quasifield, the Knuth pre-semifield and the Kantor
pre-semifield. We will  explicitly derive  formulas of divisions for
them, then three new subclasses of the $\mathcal {P}\mathcal {S}$
class can be directly obtained from Theorem \ref{mainthm}.

\subsection{The Dempwolff-M\"uller
pre-quasifield and the $\mathcal{P}\mathcal{S}_{\text{D-M}}$ class}

Assume $k$ and $m$ are odd integers with $(k,m)=1$. Let
$e=2^{m-1}-2^{k-1}-1$, $L(x)=\sum_{i=0}^{k-1}x^{2^i}$, and define a
multiplication in $\mathbb{F}_{2^m}$ as
\[x\diamond y=x^eL(xy).\]
Then $(\mathbb{F}_{2^m},+,\diamond)$ becomes a pre-quasifield,
called the Dempwolff-M\"uller pre-quasifield {\cite{dempwo}}.

To get the representation of $\diamond\frac{y}{x}$, we need the
following lemma.

\begin{lem}\label{lemDM}
Let $f(x)=x^eL(x)$ and $d(2^k-1)\equiv1\mod (2^{2m}-1)$. Then
\[f^{-1}(x)=\frac{1}{D_d(x^2)},\]
where $D_d(x)$ is the Dickson polynomial.
\end{lem}
\proof According to the proof of \cite[Theorem 3.2]{dempwo}, we know
that
\[f(x)=D_{2^k-1}\left(\frac{1}{x}\right)^{1/2}=x^{1/2}\circ D_{2^k-1}(x)\circ\left(\frac{1}{x}\right)\]
(the symbol ``$\circ$" stands for ``composition"). Then by Theorem
\ref{invdickson} we get
\[f^{-1}(x)=\left(\frac{1}{x}\right)\circ D_d(x)\circ
x^2=\frac{1}{D_d(x^2)}.\]\qedd

\begin{thm}\label{psDMclass}
Use the same notations as that of Lemma \ref{lemDM}. Then
\[\diamond\frac{y}{x}=\frac{1}{xD_d\left(\frac{y^2}{x^{2^k+1}}\right)}.\]
\end{thm}
\proof Let $F(x,y)=x\diamond y$. Then
\[F(x,y)=x^eL(xy)=\frac{1}{y^e}f(xy),\]
thus\[F_y^{-1}(x)=\frac{1}{y}f^{-1}(xy^e).\] By Lemma \ref{lemDM} we
can directly get
\[F_y^{-1}(x)=\frac{1}{y}\cdot\frac{1}{D_d(y^{2e}x^2)}=\frac{1}{yD_d\left(\frac{x^2}{y^{2^k+1}}\right)}.\]
Then the conclusion follows.\qedd

By Theorem \ref{mainthm} and Theorem \ref{psDMclass}, we obtain a
subclass of the $\mathcal {P}\mathcal {S}$ class. We call this
subclass the $\mathcal{P}\mathcal{S}_{\text{D-M}}$ class.

\subsection{The Knuth pre-semifield and the $\mathcal{P}\mathcal{S}_{\text{Knu}}$ class}
Assume $m$ is odd. For any $\beta\in\mathbb{F}_{2^m}^*$, define a
multiplication in $\mathbb{F}_{2^m}$ by
\[x\diamond y=xy+x^2\tr(\beta y)+y^2\tr(\beta x).\]
Then $(\mathbb{F}_{2^m},+,\diamond)$ becomes a pre-semifield, called
the Knuth pre-semifield {\cite{knuth}}.

To get the representation of $\diamond\frac{y}{x}$, we firstly give
the following lemma.

\begin{lem}\label{lemknu}
Let $a\in\mathbb{F}_{2^m}^*$ and
$L_a(z)=az+a^2\tr\left(\frac{1}{a}\right)z^2+\tr(z)$. Then
\[L_a^{-1}(z)=\left(1+\tr\left(\frac{1}{a}\right)\right)\frac{z}{a}+\frac{1}{a}\tr\left(\frac{z}{a}\right)+\frac{1}{a}\tr\left(\frac{1}{a}\right)C_a(z),\]
where $C_a(z)$ is the polynomial defined in Theorem \ref{invLP}.
\end{lem}
\proof When $\tr\left(\frac{1}{a}\right)=0$, we have
$L_a(z)=az+\tr(z)$, then it can be easily checked that
\[L_a^{-1}(z)=\frac{z}{a}+\frac{1}{a}\tr\left(\frac{z}{a}\right);\]
When $\tr\left(\frac{1}{a}\right)=1$, we have
$L_a(z)=az+a^2z^2+\tr(z)$, then according to Theorem \ref{invLP} we
have
\[L_a^{-1}(z)=\frac{1}{a}C_a(z)+\frac{1}{a}\tr\left(\frac{z}{a}\right).\]
Finally by Lagrange interpolation we obtain
\begin{eqnarray*}
 L_a^{-1}(z)  &=& \left(\frac{z}{a}+\frac{1}{a}\tr\left(\frac{z}{a}\right)\right)
 \left(1+\tr\left(\frac{1}{a}\right)\right)+\left(\frac{1}{a}C_a(z)+\frac{1}{a}\tr\left(\frac{z}{a}\right)\right)\tr\left(\frac{1}{a}\right) \\
   &=&\frac{z}{a}\left(1+\tr\left(\frac{1}{a}\right)\right)+\frac{1}{a}\tr\left(\frac{z}{a}\right)+\frac{1}{a}\tr\left(\frac{1}{a}\right)C_a(z).
\end{eqnarray*}\qedd

\begin{thm}\label{psknuclass}
\[\diamond\frac{y}{x}=(1+\tr(\beta x))\frac{y}{x}+x\tr\left(\beta\frac{y}{x}\right)+x\tr(\beta x)C_{\frac{1}{\beta x}}\left(\frac{y}{x^2}\right).\]
\end{thm}
\proof Let $L_a(z)$ be the polynomial defined in Lemma \ref{lemknu}
and $F(x,y)=x\diamond y$. Then
\begin{eqnarray*}
  F_y\left(\frac{x}{\beta}\right) &=& \frac{x}{\beta}y+\frac{x^2}{\beta^2} \tr(\beta y)+y^2\tr(x)\\
   &=&y^2\left[\frac{x}{\beta y}+\left(\frac{x}{\beta y}\right)^2\tr(\beta
   y)+\tr(x)\right]\\
   &=&y^2 L_{\frac{1}{\beta y}}(x).
\end{eqnarray*}
Hence by Lemma \ref{lemknu} we have
\begin{eqnarray*}
F_y^{-1}(x) &=&\frac{1}{\beta}L_{\frac{1}{\beta y}}^{-1}\left(\frac{x}{y^2}\right) \\
   &=&\frac{1}{\beta}\left[\beta y\frac{x}{y^2}(1+\tr(\beta
   y))+\beta y\tr\left(\beta y\frac{x}{y^2}\right)+\beta y\tr(\beta y)
   C_{\frac{1}{\beta y}}\left(\frac{x}{y^2}\right)
   \right]\\
   &=&\frac{x}{y}(1+\tr(\beta y))+y\tr\left(\beta\frac{x}{y}\right)+y\tr(\beta y)C_{\frac{1}{\beta
   y}}\left(\frac{x}{y^2}\right).
\end{eqnarray*}
Then the result follows from $
\diamond\frac{y}{x}=F_x^{-1}(y)$.\qedd

By Theorem \ref{mainthm} and Theorem \ref{psknuclass}, we obtain a
subclass of the $\mathcal {P}\mathcal {S}$ class. We call this
subclass the $\mathcal{P}\mathcal{S}_{\text{Knu}}$ class.

\subsection{The Kantor
pre-semifield and the $\mathcal{P}\mathcal{S}_{\text{Kan}}$ class}

Assume $m$ is odd. Define a multiplication in $\mathbb{F}_{2^m}$ by
\[x\diamond y=x^2y+\tr(xy)+x\tr(y).\]
Then $(\mathbb{F}_{2^m},+,\diamond)$ becomes a pre-semifield, called
the Kantor pre-semifield {\cite{kantor}}.

To get the representation of $\diamond\frac{y}{x}$, we firstly give
the following lemma.

\begin{lem}\label{lemkantor}
Let $L_a(z)=az^2+\tr(az)+\tr(a)z$. Then
\begin{eqnarray*}
L_a^{-1}(z)
&=&\frac{\tr(a)}{a}\left[(az)^{2^{m-1}}+\sum_{i=0}^{\frac{m-1}{2}}(az)^{2^{2i}-1}+\left(\sum_{i=0}^{\frac{m-3}{2}}a^{2^{2i}}\right)\tr(az)
\right]\\
&&+a^{2^{m-1}-1}z^{2^{m-1}}+a^{2^{m-1}-1}\tr(az).
\end{eqnarray*}
\end{lem}
\proof It is obvious that the formulation of $L_a^{-1}(z)$ is
\[L_a^{-1}(z)=\left(\frac{z}{a}\right)^{1/2}+\frac{\tr(ax)}{a^{1/2}}\]
when $\tr(a)=0$, and is
\[L_a^{-1}(z)=\frac{1}{a}\left[\sum_{i=0}^{\frac{m-1}{2}}(az)^{2^{2i}-1}+\left(\sum_{i=0}^{\frac{m-1}{2}}a^{2^{2i}}\right)\tr(az)
\right]\] when $\tr(a)=1$. Hence when $\tr(a)=0$, we have
\begin{eqnarray*}
  L_a^{-1}(L_a(z)) &=&L_a^{-1}(az^2+\tr(az))  \\
   &=&\left(z^2+\frac{\tr(az)}{a}\right)^{1/2}+\frac{\tr(a^2z^2+a\tr(az))}{a^{1/2}}\\
   &=&z+\frac{\tr(az)}{a^{1/2}}+\frac{\tr(az)+\tr(a)\tr(az)}{a^{1/2}}\\
   &=&z,
\end{eqnarray*}
and when $\tr(a)=1$, we have
\begin{eqnarray*}
  L_a^{-1}(L_a(z)) &=&L_a^{-1}(az^2+\tr(az)+z)  \\
   &=&\frac{1}{a}\left[\sum_{i=0}^{\frac{m-1}{2}}z^{2^{2i}-1}+\left(\sum_{i=0}^{\frac{m-1}{2}}a^{2^{2i}}\right)\tr(z)
\right]\circ(az)\circ\left(\frac{z^2}{a}+\frac{z}{a}+\tr(z)\right)\circ(az)\\
&=&\frac{1}{a}\left[\sum_{i=0}^{\frac{m-1}{2}}z^{2^{2i}-1}+\left(\sum_{i=0}^{\frac{m-1}{2}}a^{2^{2i}}\right)\tr(z)
\right]\circ(z^2+z+a\tr(z))\circ(az)\\
&=&\frac{1}{a}\left[\sum_{i=0}^{\frac{m-1}{2}}(z^2+z+a\tr(z))^{2^{2i}-1}+\left(\sum_{i=0}^{\frac{m-1}{2}}a^{2^{2i}}\right)\tr(z)
\right]\circ(az)\\
&=&\frac{1}{a}\left[\sum_{i=0}^{\frac{m-1}{2}}z^{2^{2i}}+\sum_{i=0}^{\frac{m-1}{2}}z^{2^{2i}-1}+\left(\sum_{i=0}^{\frac{m-1}{2}}a^{2^{2i}-1}+\sum_{i=0}^{\frac{m-1}{2}}a^{2^{2i}}\right)\tr(z)
\right]\circ(az)\\
&=&\frac{1}{a}\left(z+\tr(z)+\tr(a)\tr(z)\right)\circ(az)\\
&=&\frac{z}{a}\circ(az)\\
&=&z.
\end{eqnarray*}
 To summarize, we have $L_a^{-1}(L_a(z))=z$. (We have to remark that all operations of polynomials above are considered modulo
$(x^{2^m}+x)$.)\qedd

\begin{thm}\label{pskanclass}
\begin{eqnarray*}
\diamond\frac{y}{x}&=&\frac{\tr(x)}{x}\left[(xy)^{2^{m-1}}+\sum_{i=0}^{\frac{m-1}{2}}(xy)^{2^{2i}-1}+\left(\sum_{i=0}^{\frac{m-3}{2}}x^{2^{2i}}\right)\tr(xy)
\right]\\
&&+x^{2^{m-1}-1}y^{2^{m-1}}+x^{2^{m-1}-1}\tr(xy).
\end{eqnarray*}
\end{thm}
\proof Let $F(x,y)=x\diamond y$. Then $F(x,y)=L_y(x)$, where
$L_a(z)$ is the polynomial defined in Lemma \ref{lemkantor}. Thus
$\diamond\frac{y}{x}=L_x^{-1}(y)$ can be directly obtained from
Lemma \ref{lemkantor}.\qedd

By Theorem \ref{mainthm} and Theorem \ref{pskanclass}, we obtain a
subclass of the $\mathcal {P}\mathcal {S}$ class. We call this
subclass the $\mathcal{P}\mathcal{S}_{\text{Kan}}$ class.

\section{Concluding remarks}\label{secconclu}

In this paper, we propose a technique of constructing $\mathcal
{P}\mathcal {S}$ bent functions from finite pre-quasifield spreads.
 In fact, all $\mathcal {P}\mathcal {S}$ bent functions can be
constructed in this manner. Via this approach, we can represent all
subclasses of the  $\mathcal {P}\mathcal {S}$ class in  similar
styles with the $\mathcal {P}\mathcal {S}_{\text{ap}}$ class. To
explicitly represent these subclasses, we need to compute
compositional inverses of certain parametric permutation polynomial
over finite fields derived from multiplication operations of
pre-quasifields. Concentrated on three special classes of
pre-quasifields, we realize this goal  and obtain three new
subclasses of the $\mathcal {P}\mathcal {S}$ class, namely, the
$\mathcal{P}\mathcal{S}_{\text{D-M}}$ class, the
$\mathcal{P}\mathcal{S}_{\text{Knu}}$ class and the the
$\mathcal{P}\mathcal{S}_{\text{Kan}}$ class. They enlarge our
knowledge of functions in the $\mathcal {P}\mathcal {S}$ class more
than 30 years after it was introduced.

We have to point out that other subclasses of the $\mathcal
{P}\mathcal {S}$ class can also be obtained by our technique from
other pre-quasifields, say, the generalized Kantor pre-quasifield,
the Albert pre-quasifield, etc. Besides, we can also talk about
equivalence between these subclasses and derive the classes of their
dual functions (every bent function admits a unique dual function
{\cite{carlet2010}}). Furthermore, we can also generalize our
technique to construct bent functions over finite fields of odd
characteristic {\cite{kumar}}. All these results will be presented
in a full version of this paper.



\end{document}